\documentclass[11pt, leqno, a4paper]{amsart}
\usepackage{amsmath,amstext,amsthm,amsfonts}
\usepackage{amssymb,latexsym}
\usepackage[cp1252]{inputenc}
\usepackage[T1]{fontenc}
\usepackage{lmodern}
\usepackage{hyperref, color}
\usepackage[english]{babel}
\usepackage{mathrsfs}
\usepackage{ulem}

\theoremstyle{plain}
\newtheorem{theorem}{\textbf{Theorem}}[section]
\newtheorem{thm}[theorem]{\textbf{Theorem}}
\newtheorem{prop}[theorem]{\textbf{Proposition}}

\newtheorem{cor}[theorem]{Corollary}
\newtheorem{defi}[theorem]{Definition}

\def\boxit#1{\leavevmode\hbox{\vrule\vtop{\vbox{\kern.33333pt\hrule
    \kern1pt\hbox{\kern1pt\vbox{#1}\kern1pt}}\kern1pt\hrule}\vrule}}

\newcommand{\resp}{\emph{resp. }}
\newcommand{\ie}{\emph{i.e. }}

\newcommand{\noid}{\noindent $\diamond$~}
\newcommand{\HH}{\mathbb{H}}
\newcommand{\R}{\mathbb{R}}

\newcommand{\F}{\mathbb{F}}

\newcommand{\Z}{\mathbb{Z}}
\newcommand{\N}{\mathbb{N}}

\newcommand{\rhoh}{\widehat{\rho}}

\newcommand{\Sh}{\widehat{S}}
\newcommand{\Mh}{\widehat{M}}
\newcommand{\Dh}{\widehat{D}}
\newcommand{\Mt}{\widetilde{M}}
\newcommand{\Vh}{\widehat{V}}
\newcommand{\gh}{\hat{g}}
\newcommand{\fh}{\hat{f}}
\newcommand{\gt}{\tilde{g}}

\newcommand{\cP}{\mathcal{P}}

\newcommand{\4}{\frac{1}{4}}
\newcommand{\2}{\frac{1}{2}}
\newcommand{\lip}{\mathrm{Lip}}

\newcommand{\pf}{\par{\noindent\textbf{Proof.~}}}

\begin{document}

\normalem

\title[Spectral positivity and Riemannian coverings]{Spectral positivity
and Riemannian coverings}

\author[Pierre B\'{e}rard and Philippe Castillon]{Pierre B\'{e}rard and Philippe Castillon}

\date{01/03/2013 [To appear in Bull. London Math. Soc.]}

\begin{abstract}
Let $(M,g)$ be a complete non-compact Riemannian manifold. We
consider operators of the form $\Delta_g + V$, where $\Delta_g$ is
the non-negative Laplacian associated with the metric $g$, and $V$ a
locally integrable function. Let $\rho : (\Mh,\gh) \to (M,g)$ be a
Riemannian covering, with Laplacian $\Delta_{\gh}$ and potential
$\Vh = V \circ \rho$. If the operator $\Delta + V$ is non-negative
on $(M,g)$, then the operator $\Delta_{\gh} + \Vh$ is non-negative
on $(\Mh, \gh)$. In this note, we show that the converse statement
is true provided that $\pi_1(\Mh)$ is a co-amenable subgroup of
$\pi_1(M)$.
\end{abstract}\bigskip

\maketitle

\thispagestyle{empty}


\bigskip


\vspace{10mm}


\section{Introduction}\label{S-intro}

Let $(M^n,g)$ be a complete, connected Riemannian manifold. Denote
by $\Delta$ the non-negative Laplacian, by $\mu$ the measure
associated with the metric $g$, and consider the Schr\"odinger
operator $J=\Delta + V$, where $V$ is a locally integrable function
on $M$. Such operators appear naturally when one studies minimal (or
constant mean curvature) immersions $M^n\looparrowright N^{n+1}$.
Such hypersurfaces are critical points of a volume functional whose
second derivative is given by the quadratic form associated with the
Jacobi (stability) operator
$J=\Delta-\mathrm{Ric}_N(\nu,\nu)-|A|^2$, where $\nu$ is a unit
normal vector field along $M$ and $|A|$ the norm of the second
fundamental form of the immersion. Stable hypersurfaces are those
for which the Jacobi operator is non-negative, \ie those critical
points of the volume functional which are local minima up to second
order.

When studying stable minimal hypersurfaces, a natural question is:
``What conclusions on the Riemannian manifold $(M,g)$ can one draw
from the fact that the operator $\Delta+V$ is non-negative in the
sense of quadratic forms?'' \ie from the fact that the associated
quadratic form is non-negative on Lipschitz functions with compact
support in $M$ (or equivalently on $C^1$-functions with compact
support),
$$
0 \le \int_M \big( |df|^2 + V f^2 \big) \, d\mu \hspace{1cm} \forall
f \in \lip_0(M). \leqno{(\star)}
$$

When investigating the above question, it is often useful to pass to
the universal cover of $M$. In this paper, we study the behaviour of
non-negativity under Riemannian covering. More precisely, let $\rho
: (\Mh, \gh) \to (M,g)$ be a Riemannian covering, let $V$ be a
locally integrable function on $M$, and let $\Vh = V \circ \rho$. It
follows from \cite{FCSc80} that $\Delta_g + V \ge 0$ on $(M,g)$ (in
the sense of quadratic forms) implies that $\Delta_{\gh} + \Vh \ge
0$ on $(\Mh, \gh)$. It is a natural question to investigate the
converse statement. In \cite{MePeRo08}, Proposition~2.5, the authors
prove that the converse statement is true provided that the inverse
image $\rho^{-1}(\Omega)$ of any relatively compact subdomain
$\Omega \subset M$ has sub-exponential volume growth.\medskip

Note that this statement is not true in general and needs some
additional hypotheses on the covering, as the following simple
example shows. Consider a group $\Gamma$ acting on the hyperbolic
plane with a compact quotient $\HH^2/\Gamma$. It is well known that
the operator $\Delta - \4$ is non-negative on $\HH^2$, but it is not
the case on $\HH^2/\Gamma$ which is a compact surface.\medskip

In \cite{Bro81,Bro85}, R.~Brooks worked on a closely related
problem. Indeed, he investigated the behaviour of the infimum of the
spectrum of the Laplacian under normal Riemannian coverings. A key
assumption in his work is the amenability of the covering group. In
this paper, we prove the following result which is very much
inspired by \cite{Bro81, Bro85}.

\begin{thm}\label{T-AA}
Let $(M,g)$ be a complete Riemannian manifold of dimension $n \ge
2$. Let $\rho : (\Mh, \gh) \to (M,g)$ be a Riemannian covering. Let
$V$ be a locally integrable function on $(M,g)$, and let $\Vh =
V\circ\rho$. Assume the $\pi_1(\Mh)$ is a co-amenable sub-group of
$\pi_1(M)$. Then, the operator $\Delta_{g} + V$ is non-negative on
$(M,g)$ if and only if the operator $\Delta_{\gh} + \Vh$ is
non-negative on $(\Mh,\gh)$.
\end{thm}

As a particular case, when $\rho : (\Mh, \gh) \to (M,g)$ is a normal
covering, with covering group $G$, the assumption on the fundamental
groups of $M$ and $\Mh$ is equivalent to the amenability of the
group $G$.
\medskip

The relationship between \cite{MePeRo08}, Proposition~2.5 and
Theorem~\ref{T-AA} is that a finitely generated group with
sub-exponential growth is amenable. However, there exist amenable
groups with exponential volume growth, so that our hypothesis is
weaker than the one in \cite{MePeRo08}.\smallskip

For a complete Riemannian manifold $(M,g; V)$, equiped with a
locally integrable function $V$, consider the set
$$
I(M,g;V) = \{a \in \R ~|~ \Delta + a V \ge 0 \}.
$$
This Riemannian invariant is a closed interval which contains $0$,
see \cite{Cas06}\smallskip .

For any Riemannian covering $\rho : (\Mh, \gh) \to (M,g)$, we have
$I(M,g;V) \subset I(\Mh, \gh;\Vh)$. Theorem~\ref{T-AA} tells us that
if $\pi_1(\Mh)$ is co-amenable in $\pi_1(M)$, then $I(M,g;V) =
I(\Mh, \gh;\Vh)$: the invariant $I(M,g;V)$ does not distinguish two
manifolds which differ by a covering with amenable action. In
particular, we have the following corollary.

\begin{cor}\label{C-NC2}
Let $(M,g)$ be either $\R^2$ with is a $\Z^2$-invariant metric, or
$S^1 \times \R$ with a $\Z$-invariant metric. Let $K$ denote the
Gaussian curvature of the metric $g$. Then, $I(M,g;K) = \{0\}$ or
$I(M,g;K) = \R$. Furthermore, $I(M,g;K) = \R$ if and only if $K
\equiv 0$.
\end{cor}

\textbf{Remark}. The Jacobi operator of an isometric minimal surface
$M^2 \looparrowright \R^3$ into Euclidean $3$-space is $\Delta +
2K$, where $K$ is the Gaussian curvature of $M$. More generally
(\cite{FCSc80}, Section~3), the Jacobi operator of a minimal
immersion $M \looparrowright \Mh^3$ into a $3$-manifold with scalar
curvature $\Sh$, can be written as $\Delta + K - (\Sh + \2 |A|^2)$.
When $\Sh$ is non-negative, the stability of the surface implies
that $\Delta+K$ is non-negative. As a consequence, the study of
stable minimal surfaces led to a general spectral problem on
Riemannian surfaces, namely investigating the consequences of the
non-negativity of operators of the form $\Delta+aK$, where $K$ is
the curvature of the surface, and $a\in\R$ is a real parameter. As a
matter of fact, the papers \cite{FCSc80, Cas06, EsRo10, BerCas11}
derive topological properties of $M$ from assumptions on $I(M,g;K)$
(where $K$ is the Gaussian curvature of the metric $g$, in this
$2$-dimensional framework). \bigskip

The paper is organized as follows. In Section~\ref{S-aga}, we recall
basic facts on amenable groups and amenable group actions. For
amenable groups, we refer to \cite{Gre69}, and to \cite{Bro81},
Section 1. For amenable group actions, we refer to \cite{Ros73,
GlaMon07}. Section~\ref{S-pfs} contains the proofs.\medskip

The authors would like to thank Ivan Babenko for introducing them to
amenable group actions.


\section{Amenable group actions}\label{S-aga}

\subsection{Definitions}\label{SS-defi}

Let $X$ be a countable set, and let $G$ be a discrete group acting
on $X$. \medskip

\begin{defi}\label{D-mean}
A \emph{mean} on $X$ is a bounded linear functional $\mu :
\ell^{\infty}(X) \to \R$ such that,
$$
\mu(1)=1 \text{ ~and~ } \mu(f) \ge 0, \text{ ~whenever~} f\ge 0.
$$
Equivalently, a mean on $X$ is a map on subsets of $X$, $\mu :
\cP(X) \to [0,1]$, such that,
$$
\left\{%
\begin{array}{ll}
\mu(X)=1, \text{ ~and}\\
\mu(A\cup B) = \mu(A) + \mu(B) \text{ ~whenever~ } A\cap B = \emptyset,\\
\end{array}%
\right.
$$
\ie $\mu$ is a finitely additive probability on $X$.
\end{defi}

\begin{defi}\label{D-aag}
A $G$-action on $X$ is \emph{amenable} if there exists a
$G$-invariant mean on $X$, \ie if there exists a mean $\mu$ on $X$
such that,
$$
\mu (\gamma\!\cdot\! f) = \mu(f), ~\forall \gamma \in G, ~\forall f
\in \ell^{\infty}(X),
$$
where $(\gamma\!\cdot\! f)(x) = f(\gamma^{-1}x)$. The group $G$ is
\emph{amenable} if its action on itself by left multiplications is
amenable.
\end{defi}

Any finite group $G$ is amenable. If suffices to take $\mu(f) =
\frac{1}{\sharp (G)} \sum_{\gamma \in G} f(\gamma)$, where $\sharp
(G)$ is the number of elements of $G$. The group $\Z$ is amenable.
To see this, consider the means $\mu_n, n\in \N$, defined by
$\mu_n(f) = \frac{1}{2n+1} \sum_{k=-n}^n f(k)$ and take their
weak-$*$ limit. More generally, the group $\Z^n$ is amenable for all
$n$. The basic example of a non-amenable group is the free group on
two generators $\F_2$, and it can be shown that any group containing
a subgroup isomorphic to $\F_2$ is not amenable. \medskip

\textbf{Remarks}.\\
(i) If the group $G$ is amenable, then any action of $G$ is
amenable. Indeed, let $\nu$ be an invariant mean on $G$, and let $X$
be a set endowed with a $G$-action. Choose $x_0 \in X$. Given $f \in
\ell^{\infty}(X)$, define $\hat{f} \in \ell^{\infty}(G)$ by
$\hat{f}(\gamma)=f(\gamma\!\cdot\!x_0)$. Then, the map $\mu :
\ell^{\infty}(X) \to \R$ defined by $\mu(f) = \nu(\hat{f})$ is an
invariant mean on $X$ (we take the
mean on a $G$-orbit).\\
(ii) Note that the converse statement is not true: there exist
amenable actions by non-amenable groups. Indeed, van Douwen proved
that any finitely generated non-abelian free group admits a faithful
transitive amenable action (\cite{GlaMon07}, Introduction).

\begin{defi}\label{L-cag}
A subgroup $H < G$ is \emph{co-amenable} if the $G$-action
on $G/H$ is amenable.
\end{defi}

\textbf{Remark}. Note that a normal subgroup $H \vartriangleleft G$
is co-amenable if and only if $G/H$ is an amenable group.\medskip

When the subgroup $H$ is not normal, one can consider left cosets
$[\sigma]_L = \sigma H \in G/H$, and right cosets $[\sigma]_R =
H\sigma \in H\backslash G$. The group $G$ acts on both sets. The
$G$-action on $G/H$ is defined by $\gamma\!\cdot\![\sigma]_L = [\gamma
\sigma]_L$; the $G$-action on $H\backslash G$ is defined by
$\gamma\!\cdot\![\sigma]_R = [\sigma \gamma^{-1}]_R$. There is a natural
bijective map between the coset spaces,
$$
\left\{%
\begin{array}{l}
F: G/H \to H\backslash G,\\
F: [\sigma]_L \mapsto  [\sigma^{-1}]_R,
\end{array}%
\right.
$$
and this map is equivariant for the actions of $G$. It follows
(\cite{GlaMon07}, Lemma~2.1) that,

\begin{prop}\label{P-coa}
The $G$-action on $G/H$ is amenable if and only if the $G$-action on
$H\backslash G$ is amenable. It follows that a subgroup $H < G$ is
co-amenable if any (and hence both) of the actions of $G$ on the
coset spaces is amenable.
\end{prop}

\pf Since the map $F$ is equivariant, the push-forward of a
$G$-invariant mean on $G/H$ provides an invariant mean on
$H\backslash G$ and conversely. \qed \medskip

\subsection{F{\o}lner's condition}\label{SS-fc}

F{\o}lner's characterization of amenable groups (\cite{Gre69},
Section~3.6) has an analog for amenable actions.

\begin{thm}[F{\o}lner's condition]\label{T-fc}
A $G$-action on $X$ is amenable if and only if,
$$
\left\{%
\begin{array}{l}
\forall \eta \in (0,1), ~~\forall \gamma_1, \ldots, \gamma_n \in G,
~~\exists E \subset X  \text{ ~such that,}\\[5pt]
\sharp(E) < \infty \text{~~and~~} \eta \; \sharp(E) \le \sharp(E\cap
\gamma_i\!\cdot\!E), ~~\forall 1\le i \le n.\\
\end{array}%
\right.
$$
Equivalently, a $G$-action on $X$ is amenable if and only if,
$$
\left\{%
\begin{array}{l}
\forall \epsilon > 0, ~~\forall \gamma_1, \ldots, \gamma_n \in G,
~~\exists E \subset X  \text{ ~such that,}\\[5pt]
\sharp(E) < \infty \text{~~and~~} \sharp(E\, \triangle\,
\gamma_i\!\cdot\!E)
\le \epsilon \; \sharp(E), ~~\forall 1\le i \le n.\\
\end{array}%
\right.
$$
Here $\sharp(E)$ denotes the number of elements of $E$, and
$A\,\triangle\, B$ the symmetric difference of two sets $A, B$.
\end{thm}

\pf See \cite{Ros73}, Section~4, Theorems~4.4 and~4.9.\qed
\medskip

\textbf{Remarks}.\\
(1)~ When $G$ is countable, the second condition is equivalent to
saying that there exists a sequence $\{E_k\}_{k\in \N}$ of finite
subsets of $X$ such that,
$$
\lim_{k\to \infty} \frac{\sharp(E_k\,\triangle\,
\gamma\!\cdot\!E_k)}{\sharp(E_k)} = 0, ~~\forall \gamma \in G.
$$
This is the formulation in \cite{Ros73}.\\
(2)~ As a consequence of this criterion, one can show (\cite{Bro81},
Proposition~1) that groups with sub-exponential growth are amenable.
This fact relates our Theorem~\ref{T-AA} to Proposition~2.5 in
\cite{MePeRo08}.\medskip

For finitely generated groups, F{\o}lner's property has the following
consequence. Let $\alpha_1, \ldots, \alpha_n$ be a symmetric system
of generators for $G$. Define the Cayley graph of
the $G$-action on $X$ as follows,\\
- the vertices of the graph are the elements of $X$,\\
- $[x,y]$ is an edge of the graph if and only if there exists some
$i \in \{1, \ldots, n\}$, such that $x =
\alpha_i\!\cdot\!y$.\smallskip

For a finite set $E \subset X$, define the \emph{boundary} $\partial
E$ of $E$ (in the Cayley graph) as,
$$
\partial E = \big\{ x \in E ~|~ \exists i \in \{1, \ldots, n\},
\text{ ~such that~ } \alpha_i\!\cdot\!x \not \in E \big\}.
$$

\begin{prop}\label{P-fcg}
If the $G$-action on $X$ is amenable, then
$$
\forall \epsilon > 0, ~~\exists E \subset X, \text{ ~finite, such
that~ } ~\frac{\sharp(\partial E)}{\sharp(E)} \le \epsilon .
$$
\end{prop}

\pf Let $E$ be a finite subset of $X$. Then
$$
\partial E = \cup_{i=1}^n \{x \in E ~|~ \alpha_i\!\cdot\!x \not \in E\}.
$$
Then,
$$
\begin{array}{lll}
\sharp (\partial E) & \le & \sum_{i=1}^n \sharp (\{x \in E ~|~
\alpha_i\!\cdot\!x \not \in E\}) \\
& \le & \sum_{i=1}^n \Big( \sharp (E) - \sharp \big( \{x \in E ~|~
\alpha_i\!\cdot\!x \in E\} \big) \Big)\\
& \le & \sum_{i=1}^n \Big( \sharp (E) - \sharp (E \cap \alpha_i^{-1}
\!\cdot\!E)\Big).
\end{array}
$$

Choose $\epsilon > 0$ and $\eta \in (0,1)$ such that $n (1-\eta) \le
\epsilon $. Apply F{\o}lner's theorem: there exists a finite set $E
\subset X$ such that $\eta \; \sharp (E) \le \sharp (E \cap
\alpha_i^{-1}\!\cdot\!E)$ for $i = 1, \cdots, n$. Using the above
inequalities, we obtain for this subset,
$$
\sharp (\partial E) \le \sum_{i=1}^n \big( \sharp (E) - \eta \sharp
(E) \big) = n(1-\eta)\; \sharp (E) \le \epsilon\; \sharp (E).
$$
This proves the proposition. \qed \medskip

\textbf{Remark}. A geometric interpretation of the previous
proposition is that there is no linear isoperimetric inequality on
the Cayley graph of an amenable action by a finitely generated
group.

\section{Proofs}\label{S-pfs}

As we already mentioned in the Introduction, it follows from
\cite{FCSc80} that for any Riemannian covering $\rho : (\Mh,\gh) \to
(M,g)$, and any locally integrable function $V$ on $M$, with $\Vh =
V\circ \rho$, $\Delta_g + V \ge 0$ on $M$ implies that $\Delta_{\gh}
+ \Vh \ge 0$ on $\Mh$. We shall therefore concentrate on the
converse statement. We begin the proof by considering normal
coverings and then explain how to handle the general case.

\subsection{Normal coverings}\label{ss-pf-t-nc}

We first assume that $\pi_1(\Mh)$ is a normal subgroup of
$\pi_1(M)$. Therefore, $\rho : (\Mh, \gh) \to (M,g)$ is a normal
Riemannian covering, with amenable covering group
$G=\pi_1(M)/\pi_1(\Mh)$. In view of the introduction to this
section, we only need to prove that $\Delta_{\gh} + \Vh \ge 0$
implies that $\Delta_g + V \ge 0$.
\medskip

Let $f$ be a function in $C_0^1(M)$. We want to prove that $0 \le
\int_M |df|^2 + V f^2$. The general idea is as follows. Lift $f$ to
$\hat{f}$ on $\Mh$. This function is not compactly supported, but it
behaves like $f$ on fundamental domains. Multiply $\hat{f}$ by a
cut-off function $\xi$ which is $1$ inside some $\Omega$ to be
chosen later on. By assumption, $\int_{\Mh} |d(\xi \hat{f})|^2 + \Vh
\xi^2 \hat{f}^2$ is non-negative. We will conclude using
Proposition~\ref{P-fcg} which tells us that the effect of the
cut-off function is negligible.\medskip

\noid Fix some $\epsilon > 0$. Let $F \subset \Mh$ be a fundamental
domain for the action of the covering group $G$, and let $\beta_1,
\cdots, \beta_n$ be the elements of $G$ such that $\beta_i \!\cdot\!
F \cap F \not = \emptyset$. These elements generate $G$. As $G$ is
amenable, there exists a finite set $E \subset G$ such that $\sharp
(\partial E) \le \epsilon \, \sharp (E)$ (Proposition~\ref{P-fcg}).
\medskip

\noid Lift $\mathrm{supp}(f)$ to $F$. As it is compact, there exists
some $\alpha > 0$ such that the $\alpha$-neighborhood of $\partial
\big( F \bigcup \cup_{i=1}^n \beta_i\!\cdot\! F \big)$ does not
intersect the lift of $\mathrm{supp}(f)$. Let $\Omega = \cup_{\gamma
\in E} \, \gamma \!\cdot\! F$ and consider the cut-off function $\xi
: \Mh \to \R$ defined by
$$
\xi(x) = \left\{%
\begin{array}{l}
0, \text{ ~if~ } x \not \in \Omega, \\[5pt]
\frac{1}{\alpha}d(x, \partial \Omega), \text{ ~if~ } x \in \Omega
\text{ ~and~ } d(x, \partial \Omega) < \alpha, \\[5pt]
1, \text{ ~if~ } x \in \Omega \text{ ~and~ } d(x, \partial \Omega)
\ge \alpha.
\end{array}%
\right.
$$

The function $\xi$ satisfies the inequalities,
$$
\left\{
\begin{array}{l}
0 \le \xi \le 1 \text{ ~on~ } \Omega,\\[5pt]
|d\xi| \le \frac{1}{\alpha} \text{ ~on~ }
\gamma \!\cdot\! F, \text{ ~if~ } \gamma \in \partial E,\\[5pt]
\xi = 1 \text{ ~and~ } |d\xi| = 0  \text{ ~on~ } \gamma \!\cdot\! F,
\text{ ~if~ } \gamma
\in E\!\setminus\!\partial E.\\
\end{array}
\right.
$$

\noid Call $\fh$ the lift of $f$ to $\Mh$, and consider the test
function $\xi \,\fh$ on $\Mh$. From the assumption $\Delta_{\gh} +
\Vh \ge 0$, we have that
$$
0 \le \int_{\Mh} |d(\xi \fh)|^2 + \Vh \fh^2 \xi^2. \leqno{(a)}
$$

\noid Let $c := \sharp (E)$ and $b := \sharp (\partial E)$. Then, $b
\le \epsilon \, c$. We have
$$
|d(\fh\xi)|^2 \le \fh^2 |d\xi|^2 + 2 |\fh\xi| |d\fh| |d\xi| + \xi^2
|d\fh|^2,
$$
and it follows from the above estimates on $\xi$ that,
$$
\int_{\Mh} |d(\fh\xi)|^2 \le \frac{b}{\alpha^2} \int_M f^2 +
\frac{2b}{\alpha} \big( \int_M f^2 \big)^{1/2} \big( \int_M
|df|^2\big)^{1/2} + c \int_M |df|^2. \leqno{(b)}
$$

Consider the positive and negative parts of $\Vh$, $\Vh = \Vh_{+} -
\Vh_{-}$. Then,
$$
\begin{array}{lll}
\int_{\Mh} \Vh \fh^2 \xi^2 & = & \int_{\Mh} \Vh_{+} \fh^2 \xi^2 -
\int_{\Mh} \Vh_{-} \fh^2 \xi^2 \\[6pt]
& \le & c \int_M V_{+} f^2 - (c-b) \int_M V_{-} f^2 \\[6pt]
& \le & c \int_M V f^2 + b \int_M V_{-} f^2. \\
\end{array}\leqno{(c)}
$$

The inequalities (a)--(c) yield,
$$
0 \le \int_M |df|^2 + V f^2 + \frac{b}{c} \Big[ \frac{1}{\alpha^2}
\int_M f^2 + \frac{2}{\alpha} \big( \int_M |df|^2 \, \int_M f^2
\big)^{1/2} + \int_M V_{-}f^2 \Big].
$$

Since $0 \le b/c \le \epsilon$, letting $\epsilon$ tend to zero, we
find that $0 \le \int_M |df|^2 + V f^2$. \qed

\subsection{General coverings}\label{ss-pf-t-aa}

We first make some preparation and then explain how to make the
proof for normal coverings work in the general case.\medskip

\noid Let $\rho : (\Mh,\gh) \to (M,g)$ be a Riemannian covering. Let
$(\Mt, \gt)$ be the common universal covering for $(\Mh,\gh)$ and
$(M,g)$. Let $H = \pi_1(\Mh)$, \resp $G = \pi_1(M)$, denote the
fundamental groups of $\Mh$ and $M$ respectively. We make the
assumption that $H$ is a co-amenable subgroup of $G$.

\noid Choose a reference point $x_0 \in \Mt$. Let $\Dh$, \resp $D$,
denote the closed Dirichlet fundamental domain for the action of
$H$, \resp for the action of $G$, on $\Mt$, centered at $x_0$. As $H
< G$, we have $D \subset \Dh$. Let $\beta_1, \ldots, \beta_n$ be the
elements of $G$ such that $\beta_i D \cap D \not = \emptyset$ in
$\Mt$. These elements form a symmetric system of generators for $G$.

\noid Let $X = Gx_0 \cap \Dh$ denote the subset of points of the
orbit of $x_0$ under the action of G which belong to $\Dh$. For $x
\in X$, let $D_x$ denote the closed Dirichlet fundamental domain for
the action of $G$ on $\Mt$, centered at $x$. Then, $\Dh = \cup_{x\in
X} D_x$, and the interiors $D_x^{\circ}$ are pairwise disjoint.

\noid Let $\rhoh : \Mt \to \Mh$ be the universal covering of $\Mh$.
For $x \in X$, let $\Omega_x = \rhoh(D_x)$. Then, $\Mh = \cup_{x\in
X}\Omega_x$, and the interiors $\Omega_x^{\circ}$ are pairwise
disjoint. Since $\rho \circ \rhoh$ is the universal covering $\Mt
\to M$, the map $\rho$ is injective on each $\Omega_x^{\circ}$, and
surjective from each $\Omega_x$ onto $M$.

\noid At this point, we have tiled $\Mh$ by fundamental domains of
the covering $\rho : \Mh \to M$, indexed by the set $X$. We shall
now define an amenable action of $G$ on $X$. There is a natural map
$\Phi : X \to H\backslash G$. Indeed, given $x \in X$, there is a
unique $\gamma \in G$ such that $x = \gamma x_0$. Define $\Phi(x) =
[\gamma]_R$. We claim that the map $\Phi$ is bijective. Let $\Psi :
H\backslash G \to X$ be defined as follows. Take $H\gamma =
[\gamma]_R \in H\backslash G$. Then the set $\{\sigma x_0 ~|~ \sigma
\in [\gamma]_R\} = H(\gamma x_0)$ is the $H$-orbit of $\gamma x_0$.
Since $\Dh$ is a fundamental domain for the action of $H$, there
exists a unique $\sigma \in [\gamma]_R$ such that $\sigma x_0 \in
\Dh$. Let $\Psi([\gamma]_R) = \sigma x_0$. It is clear that $\Phi
\circ \Psi = \mathrm{Id}$ because $[\sigma]_R = [\gamma]_R$ in the
above construction. On the other-hand, $\Psi \circ \Phi =
\mathrm{Id}$ by unicity of $\sigma$. It follows that we can view the
tiling of $\Mh$ by the sets $\Omega$'s as indexed by $H\backslash
G$, and we now use the $G$-action on $H\backslash G$.

\noid We can now conclude the proof of Theorem~\ref{T-AA}.
From the above construction, $\Mh = \cup_{x\in H\backslash G}
\Omega_x$, and we have the amenable $G$-action on $H\backslash G$.

\noid By Proposition~\ref{P-fcg}, given any $\epsilon
> 0$, there exists a finite set $E \subset H\backslash G$, such
that, $\sharp (\partial E) \le \epsilon \; \sharp (E)$, where the
boundary is defined with respect to the generators $\{\beta_1,
\ldots, \beta_n\}$ defined above.

\noid We now define the set $\Omega = \cup_{x\in E} \Omega_x$ and we
can reproduce the proof of the normal covering case. \qed

\subsection{Proof of Corollary~\ref{C-NC2}}\label{ss-pf-c-nc2}

Since the groups $\Z$ and $\Z^2$ are amenable \cite{Gre69}, it
suffices to look at the quotient, \ie at the torus $T^2$. Assume
that $\Delta + a K \ge 0$ on $T^2$, for some $a \not = 0$. As in
\cite{FCSc80}, taking $u$ to be the constant function $\mathbf{1}$,
we find that $\int_{T^2} |d\mathbf{1}|^2 + a K \mathbf{1}^2 = 0$
because $\int_{T^2} K = 0$. Since $\Delta + a K \ge 0$, the function
$\mathbf{1}$ realizes the minimum of the Rayleigh quotient, so that
it is an eigenfunction associated with the eigenvalue $0$, and we
have $(\Delta + a K)\mathbf{1}=0$. This implies that $K \equiv 0$
since $a \not = 0$. As a consequence, if $I(g) \not = \{0\}$, then
$K \equiv 0$ and hence $I(g) = \R$. If $K \not \equiv 0$, then $I(g)
= \{0\}$. \qed

\bigskip

\vspace*{15mm}

\begin{footnotesize}
\noindent\begin{minipage}{0.4\textwidth}
    \begin{flushleft}
    {\normalsize Pierre B\'{e}rard}\\
    Universit\'{e} Grenoble 1\\
    Institut Fourier (\textsc{ujf-cnrs})\\
    B.P. 74\\
    38402 Saint Martin d'H\`{e}res Cedex\\
    France\\
    \verb+pierre.berard@ujf-grenoble.fr+\\
    \end{flushleft}
\end{minipage}
\hfill
\begin{minipage}{0.45\textwidth}
    \begin{flushleft}
    {\normalsize Philippe Castillon}\\
    Universit\'{e} Montpellier II\\
    D\'{e}pt des sciences math\'{e}matiques CC 51\\
    I3M (\textsc{umr 5149})\\
    34095 Montpellier Cedex 5\\
    France\\
    \verb+philippe.castillon@univ-montp2.fr+\\
    \end{flushleft}
\end{minipage}
\end{footnotesize}

\end{document}